
\input amstex
\documentstyle{amsppt}
\input bull-ppt
\keyedby{bull258e/clh}

\define\z{\zeta}
\define\lz{\Cal L_{\z}}
\define\RR{\Bbb R}
\define\CC{\Bbb C}
\define\p{\partial}

\topmatter
\cvol{26}
\cvolyear{1992}
\cmonth{Jan}
\cyear{1992}
\cvolno{1}
\cpgs{137-140}
\title Some Non-Analytic-Hypoelliptic Sums\\ 
of Squares of Vector Fields\endtitle
\shorttitle{Some nonanalytic hypoelliptic sums of squares}
\author Michael Christ\endauthor
\address  Department of Mathematics, University
of California, Los Angeles, California 90024\endaddress
\ml christ\@math.ucla.edu\endml
\date February 13, 1991 and, in revised form, June 10, 
1991\enddate
\thanks  Research supported in part by  grants from the 
National
Science Foundation\endthanks
\abstract  Certain second-order partial differential 
operators,
which are
expressed as sums of squares of real-analytic vector
fields in $\Bbb R^3$ and which are well known to be 
$C^\infty$
hypoelliptic, fail to be analytic hypoelliptic.\endabstract
 
\endtopmatter

\document

\heading 1. Introduction\endheading

A differential operator $L$ is said to be analytic 
hypoelliptic if
whenever $u$ is a distribution
such that $Lu$ is real-analytic in some open set $U$,
then $u$ is necessarily also real-analytic in $U$.
Elliptic operators with analytic coefficients are analytic
hypoelliptic, as are  certain classes of subelliptic 
operators
\cite{GS, M2, S, Ta, Tp, 
Tv1, Tv2}.    
It has been known for some time that many subelliptic
operators---whose solutions
are necessarily $C^\infty\text{---}\roman{nonetheless}$
fail to be analytic hypoelliptic; among 
the examples now known are 
\cite{BG, M1, He, PR, HH, CG}.
A substantial no-man's-land persists, in which neither 
alternative
has been proved, even in rather simple cases. 
In this note are announced negative results for
certain second-order operators. We hope that these 
will serve as models for larger classes of operators, 
rather than being
mere isolated examples.

In $\RR^3$ with coordinates $x,y,t$ set
$$
X=\frac{\p}{\p x},\qquad Y=\frac{\p}{\p 
y}-mx^{m-1}\frac{\p}{\p t},
\tag0
$$
and 
$$
L=X^2+Y^2,
$$
where $m\ge 2$ is an integer. Then $L$ is hypoelliptic in
the $C^\infty$ sense \cite{H1, K}; when $m=2$,
it is analytic hypoelliptic \cite{M2, Ta,
 Tv2}.
For $m\ge 3$ an odd integer, however, it is not analytic
hypoelliptic. This was proved for $m=3$ in \cite{He, PR},
and extended to larger $m$ in \cite{HH}, but by a method 
which
does not apply for $m$ even. 
In \cite{CG} it was found that
$\bar\partial_b\circ\bar\partial_b^*$
fails to be microlocally analytic hypoelliptic in the 
appropriate 
part of phase space,
on the CR manifold
$\{\Im(z_2)=[\Re(z_1)]^m\}$, 
where $m\ge 4$ is even. 
In appropriate coordinates for this manifold, 
$-\bar\partial_b\circ\bar\partial_b^*
=(X+iY)\circ(X-iY)=X^2+Y^2-i[X,Y]$, where $X,Y$ are as in 
(0).

\proclaim{Theorem 1}
For any even integer $m\ge 4$, $L$ is not analytic 
hypoelliptic.
\endproclaim

Despite the similarity to results just cited, this does
not follow from previous methods.
The proof is rooted in a phenomenon discovered for
$\bar\partial_b\circ\bar\partial_b^*$ 
in \cite{CG},  but that argument relied heavily upon
an explicit
formula for the Szeg\H o kernel \cite{N}, 
for which there appears to be no analogue in the present 
situation.

To place Theorem 1 in context,
consider two real
vector fields $X,Y$ in $\RR^3$ with analytic coefficients,
and suppose them to be linearly independent at each 
point. Say that
a point $a\in\RR^3$ is of type 2 if $X,Y,[X,Y]$ span the 
tangent
space to $\RR^3$ at $a$. A general result 
\cite{Tv2, Ta, M2} guarantees analytic hypoellipticity
at any point of type 2, leaving open the question of what
sort of degeneracy is permitted. 
The following conjecture has been suggested in a more 
general
form by Tr\`eves \cite{Tv2}:
$L=X^2+Y^2$ fails to be analytic hypoelliptic at $a$ if 
in any neighborhood of $a$ there exists
a real curve $\gamma$, with $\gamma'(0)\ne 0$, such that
\roster
\item"$\bullet$"
$\gamma(t)$ is not a  point of type 2 for any 
$|t|<\varepsilon$,
and
\item"$\bullet$"
$\gamma'(t)$ belongs to the span of $X(\gamma(t))$, 
$Y(\gamma(t))$
for every $|t|<\varepsilon$.
\endroster
One may hope that analytic hypoellipticity holds in all 
other cases.
In the special case of Theorem 1, 
the plane $x=0$ is foliated by a one-parameter
family of such curves $\gamma$.

More recently we have built on the analysis outlined below
to prove that analytic hypoellipticity breaks down
for $X^2+Y^2$, with $X=\p_x$ and $Y=\p_y-b'(x)\p_t$, 
whenever $b$ vanishes to
order exactly $m$ at some point, with 
$m\in\{3,4,5,\dots\}$.

\heading 2. Outline of proof\endheading

Let $\zeta,\tau$ be variables dual to $y,t$. Taking a 
partial Fourier
transform in these variables reduces the analysis of $L$
to that of a two-parameter family of ordinary 
differential operators:
$$
-\frac{d^2}{dx^2}+(\zeta-\tau mx^{m-1})^2.
$$
A simple change of variables reduces the general case
$\tau\ne 0$ to $\tau=1$, so we set
$$
\Cal L_\z = -\frac{d^2}{dx^2}+(\z-mx^{m-1})^2.
$$
It is well known \cite{H2} that in order to prove that 
$L$ is
not analytic hypoelliptic, it suffices to demonstrate
the next result (which is already known [PR, HH] for odd 
$m$).

\proclaim{Theorem 2}
Let $m\ge 3$ be an integer.
Then there exist $\z\in\CC$ and $f\in L^\infty(\RR)$, not 
identically
equal to zero, satisfying $\lz f\equiv 0$.
\endproclaim

For then, assuming that $\z$ has strictly positive 
imaginary part, one may set
$$
F(x,y,t)=\int_1^\infty e^{i\tau t+i\tau^{1/m}\z y}
f(\tau^{1/m}x)\,d\tau
$$
in the region $y>0$. Then $F\in C^\infty$, and $LF\equiv 
0$. If $f(0)\ne 0$,
one calculates readily, via a change of the contour of 
integration, that 
$$
\bigg|\frac{\p^k}{\p t^k}F(0,1,0)\bigg|\ge \delta^{k+
1}(mk)!
$$
for some $\delta>0$.
Thus $F$ is not real-analytic. If $f(0)$ does vanish, then 
$\frac{d}{dx}f(0)\ne 0$ and essentially the same 
reasoning applies
to $\frac{\p}{\p x}\frac{\p^k}{\p t^k}F$.
It is easy to check that $\lz$ has a strictly positive 
lowest eigenvalue
for each $\z\in\RR$, and that for any $\z$ satisfying the 
conclusion
of Theorem 2, $\bar\z$ does also; so the assumption above 
is legitimate.

We have only an indirect proof of the existence of 
(infinitely many)
$\zeta$ with the property desired.
Set $\gamma=-(m-1)/2$, and $\Phi_\z(x)=\z x-x^m$.
Since the coefficient of
the first-order part of $\lz$ is zero, the Wronskian of 
any two
solutions of $\lz$ is a constant function of $x$.
In the next lemma we will have two such solutions for 
each $\z$,
so their Wronskian will be a function of $\z$ alone.

\proclaim{Lemma 3}
Let $m\ge 4$ be an even integer.
For each $\z\in\CC$ there exist functions $f^+_\z$ and 
$f^-_\z$
defined on $\RR$ which satisfy 
$\lz f^\pm_\z\equiv 0$ 
and
$$
\Big|f^\pm_\z(x)-
e^{\Phi_\z(x)}
|x|^\gamma\Big|
=O(
|e^{\Phi_\z(x)}|\cdot|x|^{\gamma-1})
\quad \text{as } x\to\pm\infty,
\tag1
$$
respectively. These functions are unique, and 
depend holomorphically on $\z$. Their Wronskian, $W\!$,
satisfies
$$
|W(\z)|\le C\exp(C|\z|^{m/(m-1)}) \quad\forall\z\in\CC
\tag2
$$
for some finite $C$ and
$$
|W(\z)|\ge\delta \exp(\delta|\z|^{m/(m-1)}) 
\quad\forall\z\in\RR,
\tag3
$$
for some $\delta>0$.
\endproclaim

Now, $W$ must have at least one zero. If not, then the 
real part
of $\log W$
would be a harmonic function on $\CC^1$ with polynomial 
growth at 
infinity, hence would be
a polynomial.
By (2) and (3), its degree would have to be $m/(m-1)$. 
But for $m\ge 3$, $m/(m-1)$ is not an integer.
\footnote{Alternatively, the Hadamard product formula
guarantees that any entire function of nonintegral order 
has
infinitely many zeros.}

If $W(\z)=0$, then $f^-_\z$ is a constant multiple of $f^+
_\z$.
Hence both decay exponentially as $x\to\pm\infty$, 
therefore
certainly remain bounded. Thus $f^+_\z$ is the function 
sought.

The same reasoning can be made to apply for odd $m\ge 3$, 
with
a suitable modification of (1). 
Further argument shows that for
any $\alpha\in\RR$, the operator $X^2+Y^2+i\alpha [X,Y]$
fails to be analytic hypoelliptic. 
Related results appear in \cite{C1, C2, C4}. 

The proof of Lemma 3 is entirely elementary;
details will appear
elsewhere \cite{C3}.
The existence of solutions $f_\zeta^\pm$ with the 
prescribed
asymptotics is a special case of a standard result in the 
theory
of ordinary differential equations with irregular 
singular points
at infinity \cite{CL}.

\enddocument